%% file: main.tex
\documentclass{article}

\usepackage{amsmath}
\usepackage{amssymb}
\usepackage{enumerate}
\usepackage{amsthm}
\usepackage{listings}
\usepackage[all]{xy}
\usepackage{longtable}
\usepackage{graphicx}
\usepackage[backend=bibtex]{biblatex}
\theoremstyle{definition}
\newtheorem{theorem}{Theorem}[section]
\newtheorem{remark}[theorem]{Remark}

\newtheorem{proposition}[theorem]{Proposition}
\newtheorem{lemma}[theorem]{Lemma}
\newtheorem{example}[theorem]{Example}%
\newtheorem{corollary}[theorem]{Corollary}
\newtheorem{construction}[theorem]{Construction}%

\newcommand{\Mat}{\operatorname{Mat}}

\newcommand{\C}{\mathbb{C}}
\newcommand{\CC}{\mathbb{C}}
\newcommand{\PP}{\mathbb{P}}
\newcommand{\sheafhom}{\mathcal{H}\kern -.5pt om}
\newcommand{\coker}{\operatorname{coker}}

\title{On the Vanishing Topology of Isolated Cohen-Macaulay 
       Codimension $2$ Singularities}
\author{Anne Fr\"uhbis--Kr\"uger, Matthias Zach\\
        Institut f. Alg. Geometrie, Leibniz Universt\"at Hannover, Germany}

\begin{document}

\maketitle

\begin{abstract}
The topology of isolated complete intersections is well studied, but beyond 
this class not much is known. Isolated Cohen-Macaulay codimension 2 
singularities share many common features with isolated complete intersection 
singularities, but they also exhibit some striking new behaviour. One such 
instance was observed by Damon and Pike \cite{DP2} in their study of the 
vanishing topology and Euler characteristic, where they took this class 
of singularities as examples. 
In this article, we explore the background and geometrical meaning of 
their findings by determining the Betti numbers explicitly and explain 
the new phenomena. An important tool here is the Tjurina modification 
relating a Cohen-Macaulay codimension 2 singularity to a finite number 
of complete intersection singularities.
\end{abstract}

\input Intro.tex

\input chapter2.tex

\input chapter3.tex

\input chapter4.tex

\newpage
\input chapter5.tex

\printbibliography

\end{document}

%% file: Intro.tex
\section{Introduction}
Isolated hypersurface singularities and, a bit more generally, isolated 
complete intersection singularities have been a central focus of 
singularity theory ever since the famous A-D-E list of Arnold \cite{Arn} 
classifying the simple hypersurface singularities. Classification questions
as well as topological and analytic properties of these singularities have 
been studied intensively over the past decades (e.g. \cite{Giu},
\cite{GS}); many of the properties can be expressed in terms of 
invariants of the singularities and relations among these. The most famous 
ones are 
the Milnor number $\mu$ on the topological side 
and the Tjurina number $\tau$, i.e. the dimension of the $T^1$, 
related to the first order deformations. 
Following Damon and Pike we will view $\mu$ as the ``defect of the 
Euler characteristic''. For an isolated singularity with unique Milnor
fiber this is the difference between the topological Euler characteristic
of any smooth fibre and the central fiber in a deformation.
For quasihomogeneous ICIS, equality of the two numbers holds (see \cite{Gr}); 
for any general ICIS we still have $\mu \geq \tau$ (see \cite{LSt}).\\

As soon as we pass beyond ICIS, however, only a few results are known. 
The easiest non-ICIS case is the case of isolated Cohen-Macaulay codimension 2 
singularities where the Hilbert-Burch theorem allows a description by means
of the presentation matrix of the vanishing ideal. This case has recently 
come into focus of ongoing research, starting with the classification
simple singularities in this case, first for space curves in \cite{FK1}
and later for arbitrary dimension in \cite{FN}. This, in turn led to
further study of the properties of these singularities as in 
\cite{dSR}, \cite{NOT} and most recently in \cite{GaR} and \cite{CS}.
But up to now the properties of the Milnor fibre of such singularities 
are far from being explored in all details.\\

Damon and Pike studied the vanishing topology of a generalization of 
the Milnor fibre, the so-called singular Milnor fibre in \cite{DP2}, in 
a much broader context. They used the simple isolated Cohen-Macaulay
codimension 2 singularities from the list by Fr\"uhbis-Kr\"uger and 
Neumer \cite{FN} as examples to illustrate their methods. For the 
surface case, their approach (as well as independently a different 
approach introduced by da Silva Pereira and Ruas in \cite{dSR}) provided 
direct computations for the Milnor number, which is the second (and only 
non-vanishing) Betti number of the Milnor fiber in this dimension. \\

Moving one dimension higher to isolated Cohen-Macaulay codimension 2  
singularities of dimension 3, the situation becomes more delicate. Some
facts are still known: the existence of a smoothing and the 
vanishing of the first Betti number have been shown by Greuel and 
Steenbrink in \cite{GS}, whereas the vanishing of homology in degrees 
bigger than the complex dimension of the underlying variety 
is a well-known consequence of the Lefschetz hyperplane theorem
(see \cite{Mil63}). The methods of Damon and Pike then allowed the 
computation of the difference $b_3-b_2$ of the two remaining Betti-numbers,
but not of each of the two separately. This is the defect of the Euler
characteristic which we understand as a generalization of the Milnor number 
$\mu$ as mentioned earlier. Nevertheless their results provided
striking evidence for $b_2$ to be nonzero in some of the families of simple 
threefold singularities from \cite{FN}, as the computed value of the 
difference was negative, but never smaller than $-1$ for these. \\

In this article, we extend a technique which was previously only used 
for surfaces, the Tjurina modification (see e.g. \cite{Tj} or \cite{DvS}),
applying it to Cohen-Macaulay codimension 2 singularities in general.
It allows us to relate the given singularity or family of singularities
to a local complete intersection scheme or a 
family of local complete intersections factoring through a given deformation. 
Using this tool, we are then able to 
explain the observation of Damon and Pike, explicitly compute that 
Betti numbers $b_2$ and $b_3$ for all simple isolated Cohen-Macauly 
codimension 2 threefold singularities and even state a large class of such 
singularities (including the simple ones) for which $b_2$ has to have the 
value $1$.  \\

In section 2 we briefly recall those 
of the known results about isolated Cohen-Macaulay codimension 2 
singularities which will be needed later on. In the section 3 , 
we consider the notion of a Tjurina modification in detail and extend 
it suitably to higher dimensions, larger matrices and families of 
singularities. In these two sections, we also give very 
explicit descriptions of the objects and statements to allow the use of 
the results in algorithmic and experimental approaches to problems of similar 
flavour. We then recall important facts about the Milnor fibre 
and prove the main results in section 4. The last section contains
the application of the results to explicit examples.\\
 
We would like to thank Terence Gaffney, Wolfgang Ebeling, Wim Veys, 
Slawomir Rams, Victor Gonzalez Alonso, Miguel Marco and Jesse 
Kass for fruitful exchange of ideas on the topics related to this article.
This work is partially supported by funds of the research project 'Experimental
methods in Computer Algebra' of the NTH.

%% file: chapter2.tex
\section{Basic Facts on isolated Cohen-Macaulay codimension 2 singularities}

\label{Basics} 
Isolated Cohen-Macaulay codimension 2 singularities (abbreviated by
ICMC2 in the following) provide 
the most accessible setting for non-complete-intersection singularities. 
Contact equivalence, semi-universal deformation and simple objects are
known in this case, see \cite{FK1} and \cite{FN}. For the list of simple
objects in the dimension 3 case, see table \ref{tab:ResultsThreefolds}.  
In this section, 
we will briefly recall some of the basic facts for reader's convenience. \\

Using the Hilbert-Burch theorem, all Cohen-Macaulay germs of 
codimension 2 can be expressed as the maximal minors of 
$(t+1)\times t$-matrices $M$ and vice versa. In the same way, flat 
deformations can be represented by perturbations of the matrix $M$ and 
any perturbation gives rise to a flat deformation (cf. Burch \cite{Bur}, 
Schaps \cite{Sch}). The minimal matrix size $t$ is called the 
Cohen-Macaulay type of the singularity. 
 
Classification up to contact-equivalence means that two singularities
are considered equivalent, if their germs are isomorphic. The action of
the contact-group translates directly to the application of coordinate
changes and row and column operations on $M$. A singularity is called
simple, if it can only deform into finitely many different equivalence
classes (types) of singularities.\\

For a more consistent notation, we prefer to describe the Cohen-Macaulay
codimension 2 singularities by their presentation matrix instead of the
vanishing ideal. This requires a reformulation of $T^1_{X,0}$ in terms of 
the presentation matrix:

\begin{lemma}[\cite{FK1}] \label{lem:structT1}
$T^1_{X,0}$ is given by
$$T^1_{X,0} \cong \Mat(t+1,t;{\mathbb C}\{x_1,\dots,x_n\}) / (J_M + \operatorname{Im}(g))$$
where $J_M$ is the submodule generated by the matrices of the form
$$ \begin{pmatrix}
  \frac{\partial M_{11}}{\partial x_j} & \dots &
  \frac{\partial M_{1t}}{\partial x_j} \cr
  \vdots & & \vdots \cr
  \frac{\partial M_{(t+1)1}}{\partial x_j} & \dots &
  \frac{\partial M_{(t+1)t}}{\partial x_j} \end{pmatrix}\;\;\; \forall 1 \leq j
  \leq m$$
and $g$ is the map 
\begin{multline*}
\Mat(t+1,t+1;\C\{x_1,\ldots,x_m\})\oplus \Mat(t,t;\C\{x_1,\ldots,x_m\})
\\ \overset{g}{\rightarrow}  \Mat(t+1,t;\C\{x_1,\ldots,x_m\})
\end{multline*}
mapping $(A,B)  \mapsto  AM+MB$.
\end{lemma}
 
It is a well-known fact that $T^2_{X,0}=0$ for Cohen-Macaulay codimension 2
singularities, i.e. that there are no obstructions to lifting first order
deformations. As the Cohen-Macauly codimension 2 singularities,
we are considering, are isolated,  $T^1_{X,0}$ is of finite dimension 
$\dim_{\mathbb C} T^1_{X,0} = \tau$. Hence the base of the semiuniversal
deformation of $(X,0)$ is $\mathbb C^\tau$ and the total space is given by the 
minors of the matrix
$$M_{su}= M + \sum_{i=1}^{\tau} s_i m_i \in 
   \Mat(t+1,t;{\mathbb C}[s_1,\dots,s_{\tau}]\{\underline{x}\})$$
where the $s_i$ are the coordinates of $\mathbb C^\tau$ and 
$\{m_1,\dots,m_\tau\}$ is a ${\mathbb C}$-basis of $T^1(X,0)$ in the
matrix notation of Lemma 1. \\

The above description of $T^1_{X,0}$ in terms of the presentation matrix, 
the non-existence of obstructions and the explicit description of the 
semiuniversal deformation are the main reasons why this class of
singularities is a natural choice for a first step beyond isolated complete 
intersection singularities: 
To study their deformations we can follow the main ideas used in the 
complete intersection case.
This led to the 
complete classification of simple isolated Cohen-Macaulay codimension 2 
singularities found in \cite{FN} 
which lists nearly 30 series and more than 20 exceptional 
cases in dimensions $0 \leq \dim(X,0) \leq 4$. 

On the other hand, there are certain structural properties of our 
singularities which do not coincide with the complete intersection case
and are based on the fact that the ring of $(X,0)$ is a determinantal 
ring (see e.g. \cite{BV} 
for a textbook on determinantal rings). A first occurrence of this
situation is linked to the following well-known fact:

\begin{lemma} \label{lem:genericMatrix} 
For $k \leq l$,
let $G \in \Mat(l,k;\C\{y_{1,1},\dots,y_{l,k}\})$ be the matrix with
entries $g_{i,j}=y_{i,j}$. 
Denote by $V_r \subset \mathbb C^{l\cdot k} = \Mat(l,k;\mathbb C)$ 
the variety of matrices with 
rank $\leq k-r$. 
These form a chain
\[ \mathbb C^{l \cdot k} = V_0 \supset V_1 \supset \cdots V_{k} = 
    \left\{ 0 \right\},
\]
where the variety $V_r$ is defined by the ideal generated by all 
$k-r+1$-minors of $G$.
For $k > r > 0$ the singular locus of $V_r$ is precisely $V_{r+1}$.
\end{lemma} 

In fact, one implication of the last statement is obvious: 
the entries of the 
jacobian matrix of the ideal of $r$-minors of $G$ are 
$\mathbb C\{\underline y\}$-linear 
combinations of the $(r-1)$-minors as we can easily check by direct 
compuation.

Now suppose an ICMC2 singularity $(X,0)\subset (\mathbb C^n,0)$ is given by a matrix 
$M \in \Mat(t+1,t; \mathbb C\{\underline x\}$. We can regard this matrix
as a map 
\[ 
    M : (\mathbb C^n,0) \to \Mat(t+1,t;\mathbb C) \cong \mathbb C^{t\cdot (t+1)},
\]
which by abuse of notation we also denote by $M$. 
Then the singularity $(X,0)$ appears as the preimage 
$M^{-1}(V_1)$.

\begin{corollary}
Let $M \in \Mat(t+1,t;\C\{x_1,\dots,x_n\})$ a presentation matrix of an 
isolated Cohen-Macaulay codimension 2 singularity $(X,0)$ of Cohen-Macaulay
type $t$. Then the 
locus defined by the $(t-1)$-minors of $M$ is either the origin or 
empty. 
\label{cor:MaxMinorsOrigin}
\end{corollary}

\begin{proof}
As $(X,0)$ is a germ of an isolated singularity, the
singular locus is the origin. 
Regarding $M$ as a map to $\mathbb C^{t\cdot (t+1)}$ gives a ring 
homomorphism
\begin{eqnarray*}
M^* : \C \{y_{i,j} \mid 1 \leq i \leq t+1, 1 \leq j \leq t\} &
       \longrightarrow & \C \{x_1,\dots,x_n\}\\
y_{i,j} & \longmapsto & m_{i,j}.
\end{eqnarray*}
Here $m_{i,j}\in \mathbb C\{x_1,\dots,x_n\}$ denotes the entry of the matrix 
$M$ in the $i$-th row and $j$-th column.
Let $\langle \delta_1,\dots,\delta_{t+1}\rangle$ be the 
ideal generated by the $t+1$ maximal minors 
$\delta_i \in \mathbb C\{\underline y\}$ of 
the matrix $G$ from the previous lemma, i.e. the vanishing ideal of
$V_1 \subset \mathbb C^{t \cdot (t+1)}$.
Then the 
ideal $I \subset \C\{\underline x\}$ defining $X$ is given by 
\[ f_i(\underline x) = \delta_i(M(\underline x)).\]
Hence the jacobian matrix factors
\[ \operatorname{J}_f(\underline x) = 
    \operatorname{J}_g (M(\underline x)) \cdot
    \operatorname{J}_M(\underline x).
\]
Now by the preceeding lemma the matrix 
$\operatorname{J}_g (M(\underline x))$ has entries contained in the 
ideal of $(t-1)$-minors of $M$ and hence does 
$\operatorname{J}_f(\underline x)$ and all ideals of minors thereof. 
The statement follows immediately from this inclusion of ideals.
\end{proof}

\begin{remark}
    Of course, Corollary (\ref{cor:MaxMinorsOrigin}) can also be proved
    more directly. However, that does not illustrate the point in question.
    This is the more elegant argument: 
    The presentation matrix $M$ of the isolated 
    singularity $(X,0)$ describes the relations of the generators 
    of the conormal module $I/I^2$. Thus the ideal of $t-1$-minors
    is the second fitting ideal 
    $\operatorname{Fitt}_2(I/I^2)$ of the conormal module. Its 
    vanishing locus is the set of primes where $I/I^2$ cannot be 
    generated by $2$ elements. But $X$ is of codimension $2$ 
    and $I/I^2$ is locally free on the smooth part. 
    Hence 
    \[V(\operatorname{Fitt}_2(I/I^2)) \subset \operatorname{Sing}(X)
        = \{0\}.
    \]
\end{remark}

\begin{remark} \label{rmk:KV}
Regarding the matrix $M$ of an ICMC2 singularity as a map 
to the space of matrices as used in Corollary (\ref{cor:MaxMinorsOrigin})
provides a different perspective 
to those of the properties of our singularities which originate from 
the determinantal structure: Cohen-Macaulay codimension 2 singularities are 
among the classes of singularities for which deformations of the space 
germs coincide with deformations of $M$ as a map
(see \cite{Buc}, chapter 4 and 5). 
The appropriate notion of equivalence in this context is 
${\mathcal K}_V$-equivalence (see \cite{Da1}, \cite{Da2}),  
but we will not need this notion in our considerations.
\end{remark}

In this article, the interplay of both aspects, i.e.
of the similarities to the ICIS case and of the structural properties of 
determinantal singularities, will be essential to studying the topology of 
the singularities in question. More precisely, we shall even see 
contributions of both kinds in the topology.

%% file: chapter3.tex
\section{Tjurina modifications revisited}

Central to our considerations will be a not so widely known tool that was 
developed by G. Tjurina in \cite{Tj} for her study of rational triple point 
singularities. After finding that such surface singularities can be realized
by a system of $3$ equations, which we can easily recognize as $2$-minors of 
a $2 \times 3$ matrix, she considers the map to ${\mathbb P^1}$ which maps 
each point of the rank-1-locus of the matrix to the corresponding (non-zero) 
column vector of the matrix. Resolving the locus of indeterminacy of this 
map (i.e. the rank-0-locus) then provides her with a local complete 
intersection which only possesses rational double point singularities.\\

This construction has later also been used in the thesis of D. van Straten 
\cite{DvS}, where its name was coined, and in a few other articles. However, 
it has -- to our knowledge -- never been applied beyond the case of surface 
singularities of Cohen-Macaulay-type $t=2$. As we shall apply it to the
3-dimensional case and as we do not want to restrict our methods to the
case of Cohen-Macaulay-type $t=2$, we will generalize Tjurina's construction
here.\\

\begin{construction}{(Tjurina modification for 
        generic determinantal varieties)}
        \label{cst:tjurinaModificationDet}

    Consider the varieties $V_r \subset \mathbb C^{l \cdot k},\; k\leq l$,  
    as in Lemma \ref{lem:genericMatrix}. 
    For a general point, i.e. a general $l\times k$-matrix $A \in V_r$,
    the row vectors of $A$ span a
    $k-r$ dimensional hyperplane $P_A \subset \mathbb C^k$. 
    This determines a rational map 
    to the Grassmannian of $(k-r)$-planes in $k$-space.
    \begin{eqnarray*} 
        P: V_r & \dashrightarrow & \operatorname{Grass}(k-r,k)\\
        A & \mapsto & P_A.
    \end{eqnarray*}
    Clearly $P$ is defined on the open set $V_r\setminus V_{r+1}$.
    Regarding $\operatorname{Grass}(k-r,k) \subset 
    \mathbb P( \bigwedge^{k-r} \mathbb C^k )$ as 
    a subvariety of projective space it becomes clear that 
    $P$ can always be expressed in terms of $k-r$-minors of $A$.
    As a projective variety the Grassmannian is complete and
    we can blow up the rational map $P$ to obtain
    \[ 
        W_r := \overline{ \Gamma_{P}(V_r\setminus V_{r+1})} 
        \subset \mathbb C^{l\cdot k} \times \operatorname{Grass}(k-r,k)
    \]
    as the closure of the graph of $P$ restricted to 
    $V_r \setminus V_{r+1}$ together with the canonical 
    projection $\pi$ and the prolongation $\hat P$

    \[
        \begin{xy}
            \xymatrix{& W_r \ar[d]^{\pi} \ar[dr]^{\hat P} &  \\
            & V_r \ar@{-->}[r]^{P} & \;\;\;\operatorname{Grass}(k-r,k)}
        \end{xy}
    \]
    \noindent
    In particular $\pi$ is a resolution of the singularities of $V_r$.
\end{construction}

\begin{remark} 
    For calculating $W_r$ explicitly, we cover the projective variety
    $\operatorname{Grass}(k-r,k) \subset {\mathbb P}(\bigwedge^{k-r} \mathbb C^k)$
    by the standard affine charts.
    Similarly to writing a point $p\in \mathbb P^n$ in projective $n$-space
    as $p = (s_0:\dots:s_n)$ in projective coordinates and thus also indicating
    the line 
    $L(p) = \operatorname{span}((s_0,\dots,s_n)^T) \subset \mathbb C^{n+1}$ 
    sitting over
    $p$ in the tautological bundle, we write a point 
    $z \in \operatorname{Grass}(k-r,k)$ as a $(k-r) \times k$-matrix $B$. 
    The standard cover is indexed by subsets
    $\alpha \subset \{ 1,\dots,k \}$ of cardinality $\# \alpha = k-r$. 
    Analogous to normalizing the projective coordinates of 
    a point $p=(s_0:\dots:s_i:\dots:s_n)$ in the $i$-th chart of 
    projective space to 
    \[ 
        p = \left( \frac{s_0}{s_i}:\dots:1:\dots:\frac{s_n}{s_i} \right)
        = \left( s_0^{(i)}: \dots: 1: \dots :s_n^{(i)} \right), 
    \]
    we require the maximal square submatrix of $B$ indexed by $\alpha$ 
    to be the unit matrix. Thus we obtain affine coordinates 
    $(z^{(\alpha)}_{i,j})_{i,j}$.
    For example if $\alpha = \{ 1,\dots,k-r\}$ we write a point 
    $z \in U_\alpha \subset \operatorname{Grass}(k-r,k)$ as 
    \[
        B_\alpha(z) = 
        \begin{pmatrix}
            1 & 0 & \cdots & 0 & z^{(\alpha)}_{1,k-r+1} & \cdots & z^{(\alpha)}_{1,k} \\
            0 & 1 & \ddots & \vdots & z^{(\alpha)}_{2,k-r+1} & \cdots & z^{(\alpha)}_{2,k} \\
            \vdots & \ddots & \ddots & 0 & \vdots & & \vdots \\
            0 & \cdots & 0 & 1 & z^{(\alpha)}_{k-r,k-r+1} & \cdots & z^{(\alpha)}_{k-r,k} \\
        \end{pmatrix}
    \]
    The subspace $L(z) \subset \mathbb C^k$ sitting over $z$ is now given by the 
    span of the rows of the above matrix. 
    Given a generic $l\times k$ matrix $G\in V_r \subset \mathbb C^{l\cdot k}$,
    requiring the span of the rows of $G$ to be contained in $L(z)$ therefore 
    amounts to asking for the $k-r+1$-minors of the matrix 
    $$\left ( \begin{array}{c} B_{\alpha} \cr \hline \cr G \end{array} \right )$$
    to vanish.  
    In fact the variety 
    $W_r \subset \mathbb C^{l\cdot k} \times \operatorname{Grass}(k-r,k)$
    which is locally defined by these minors is already the strict transform 
    of $V_r$ under the blowup of $P$ in our construction.
\end{remark}

In the setting of ICMC2 singularities we will only be concerned with 
$(t+1) \times t$-matrices of rank $t-1$. A $t-1$-dimensional subspace
$L$ 
in $\mathbb C^t$ is uniquely determined by the class of a normal
vector $[\vec n_L ] \in \mathbb P^{t-1} \cong \operatorname{Grass}(t-1,t)$.
The identification of $\operatorname{Grass}(t-1,t)$ with $\mathbb P^{t-1}$
is given on the standard cover by identifying a point\footnote{The 
    non-standard numbering
    was chosen for consistency with the numbering of matrix entries.}
\[(s_1^{(i)}:\dots:s^{(i)}_{i-1}: 1 : s^{(i)}_{i+1} : \dots :s^{(i)}_t) \] 
in the chart $ \{s_i \neq 0\}$ in $\mathbb P^{t-1}$
with the matrix 
\begin{equation}
    B_i(s) = 
    \begin{pmatrix}
        1 & 0 & \cdots & 0 & -s^{(i)}_1 & 0 & \cdots & \cdots & 0 \\
        0 & 1 & \ddots & \vdots & \vdots & \vdots & & & \vdots \\
        \vdots & \ddots & \ddots & 0 & \vdots & \vdots & & & \vdots \\ 
        0 & \cdots & 0 & 1 & -s^{(i)}_{i-1} & 0 & \cdots & \cdots & 0 \\
        0 & \cdots & \cdots & 0 & -s^{(i)}_{i+1} & 1 & 0 & \cdots & 0 \\
        \vdots & & & \vdots & \vdots & 0 & 1 & \ddots & \vdots \\
        \vdots & & & \vdots & \vdots & \vdots & \ddots & \ddots & 0 \\
        0 & \cdots & \cdots & 0 & -s^{(i)}_t & 0 & \cdots & 0 & 1
    \end{pmatrix}
    \label{eqn:matrixChartCodimensionOne}
\end{equation}
The equations of the Tjurina transform $W_1$ of 
$V_1 \subset \mathbb C^{t\cdot(t+1)}$ therefore take a particular simple
form: For a point $s = (s_1:\dots:s_t) \in \mathbb P^{t-1}$ we just 
require the vector $\vec s = (s_1,\dots,s_t)^T$ to be perpendicular to the 
columns of a matrix $G \in V_1$.

\begin{corollary} 
    Let $G=(y_{i,j})_{i,j} \subset \C\{\underline{y}\}$ be the generic
    $(t+1) \times t$ matrix and $(s_1:\dots:s_t)$ the projective coordinates
    of $\mathbb P^{t-1} = \operatorname{Grass}(t-1,t)$.
    The Tjurina transform
    $W_1 \subset \mathbb C^{(t+1)t} \times \mathbb P^{t-1}$ of $V_1$ 
    is the zero locus of the equations
    \begin{equation}
        \begin{pmatrix}
            y_{1,1} & \cdots & y_{1,t} \\ 
            \vdots &  & \vdots \\
            y_{t+1,1} & \cdots & y_{t+1,t}
        \end{pmatrix}
        \cdot
        \begin{pmatrix}
            s_1 \\ \vdots \\ s_t
        \end{pmatrix}
        = 
        \begin{pmatrix}
            0 \\ \vdots \\ 0
        \end{pmatrix}.
        \label{eqn:EquationsTjurinaModCodOne}
    \end{equation}
\end{corollary}

\begin{proof}
    Computing each chart with the matrix $B_i$ as in 
    (\ref{eqn:matrixChartCodimensionOne}) 
    gives the local equations. But these can be easily recognized as the 
    dehomogenization of the equations given by 
    (\ref{eqn:EquationsTjurinaModCodOne}).
\end{proof}

From now on we will denote a representative of an isolated singularity 
by $X_0$, i.e. with an additional index $0$. 
This is due to consideration of deformations in the sequel,
where the singularity $(X_0,0)$ is embedded as the 
special fiber in a total space 
$(X,0) \overset{\varepsilon}{\longrightarrow} (B,0)$ 
fibered over some base $(B,0)$ by deformation parameters
$\varepsilon$. Consequently for a choice of representatives
the fiber over any nonzero $\varepsilon \in B$ is denoted
by $X_\varepsilon$.

\begin{construction}{(Tjurina modification for ICMC2 singularities)}
    \label{cst:tjurinaModificationICMC2}
As pointed out in remark \ref{rmk:KV} an ICMC2 
singularity $(X_0,0) \subseteq \mathbb ({\mathbb C}^m,0)$
can be studied by means of a polynomial map $M: U \longrightarrow
{\mathbb C}^{(t+1)t}$ for a chosen representative 
$X_0 \subset U$ of $(X_0,0)$ in a neighborhood $U$ of the origin \footnote{By 
abuse of notation we will also refer to such a representative as ICMC2.}.
Concatenation with $P : V_1 \dashrightarrow \mathbb P^{t-1}$
gives a rational map 
\[ P \circ M : X_0 \dashrightarrow \mathbb P^{t-1}. \]
Because $P \circ M$ is expressed in the projective coordinates of 
$\mathbb P^{t-1}$ in terms of $(t-1)$ minors of $M$,
it is well defined outside the singular locus of $X$ by Corollary 
\ref{cor:MaxMinorsOrigin}.

Now we define the Tjurina modification $Y_0$ to be the fiber product
$X_0 \times_{V_1} W_1$ in the following diagram:
\begin{equation}
    \begin{xy}
        \xymatrix{
            X_0\times_{V_1} W_1 \ar[r]^{\hat M} \ar[d]_\pi &
            W_1 \ar[d]_\rho \ar[dr]^{\hat P} &
            \\
            X_0 \ar[r]^M &
            V_1 \ar@{-->}[r]^{P}&
            {\mathbb P}^r\\
        }
    \end{xy}
    \label{eqn:DiagramTjurinaMod}
\end{equation}
On the level of equations this means nothing but regarding $M$
as a matrix with polynomial entries and requiring the equations of the system 
\[ M \cdot
    \begin{pmatrix}
        s_1 \\ \vdots \\ s_t
    \end{pmatrix}
    = 0
\]
to hold in $U \times \mathbb P^{t-1}$. 
Clearly outside the singular locus $\{0\}$ the map 
$\pi : Y_0 \to X_0$ is an isomorphism, while the origin itself is substituted 
with the whole Grassmannian $\mathbb P^{t-1}$. 
\end{construction}

\begin{example} \label{34example}
Let us consider the (non-simple) ICMC2 singularity
$(X_0,0) \subset ({\mathbb C}^5,0)$ given by the $3$-minors of the matrix 
\[
    M = 
    \begin{pmatrix}
        x & y-v & y+z \cr
        y & z-v & x+u \cr
        z & 0 & x-u \cr
        0 & u & v
    \end{pmatrix}.
\]
Let $(s_1:s_2:s_3)$ be the projective coordinates 
of ${\mathbb P}^2 = \operatorname{Grass}(2,3)$.
Then we obtain $Y_0 \subset \mathbb C^5 \times \mathbb P^2$ as 
the zero locus of the equations
\begin{equation}
    \begin{pmatrix}
        x & y-v & y+z \cr
        y & z-v & x+u \cr
        z & 0 & x-u \cr
        0 & u & v
    \end{pmatrix} 
    \cdot
    \begin{pmatrix}
        s_1 \\ s_2 \\ s_3
    \end{pmatrix}
    = 0.
    \label{eqn:Example4by3}
\end{equation}
The Tjurina transform $Y_0$ is still singular at 10 distinct points in 
$\{0\} \times {\mathbb P}^2 \subset \mathbb C^5 \times \mathbb P^2$.
But there we only find 3-dimensional $A_1$ singularities embedded in higher 
dimensional space. Thus the situation became much simpler.
Consider e.g. the singularity at the point $p = (0,(1:0:0))$ in the chart 
$s_1\neq 0$: 
The first three lines of the system 
(\ref{eqn:Example4by3}) define a smooth variety $H$ of dimension
$4$ around $p$. Inside $H$ the equation 
\[ s^{(1)}_2 \cdot u + s^{(1)}_3\cdot v = 0 \]
in the last line provides the $A_1$ singularity.
\end{example}

Any deformation of an ICMC2 singularity $X_0$ is described by a
perturbation of the entries of the matrix $M$ defining $X_0$. 
As the process of 
Tjurina modification is based on the determinantal structure, it can be
applied to all fibers of such a family simultaneously. Therefore it makes 
sense to ask whether (or in which situations) a Tjurina modification is
well-behaved within the family. 

\begin{construction}{(Tjurina modification in family)}
    Let 
    $X_0 \hookrightarrow X \overset{\varepsilon}{\longrightarrow} \mathbb C$ 
    be a deformation of an ICMC2 singularity $X_0 \subset \mathbb C^n$ 
    of Cohen-Macaulay type $t$ described by a matrix
    $M(\underline x, \varepsilon) 
    \in \Mat(t+1,t;\mathbb C\{\underline x, \varepsilon\})$. 
    The Tjurina modification in family for this deformation is the 
    result of applying the Tjurina modification to the total space 
    $X\overset{\varepsilon}{\longrightarrow} \mathbb C$ which leads to a 
    diagram extending diagram (\ref{eqn:DiagramTjurinaMod}) above:
    \begin{equation}
        \begin{xy}
            \xymatrix{
                X_0 \times_{V_1} W_1 \ar@{^{(}->}[r] \ar[d]_{\pi_0} &
                X\times_{V_1} W_1 \ar[r]^{\hat M} \ar[d]_\pi &
                W_1 \ar[d]_\rho \ar[dr]^{\hat P} &
                \\
                X_0 \ar@{^{(}->}[r] \ar[d] & 
                X \ar[r]^M \ar[d]_\varepsilon &
                V_1 \ar@{-->}[r]^{P}&
                {\mathbb P}^r
                \\
                \{0\} \ar[r] & 
                \mathbb C
            }
        \end{xy}
        \label{eqn:DiagramTjurinaModInFamily}
    \end{equation}
    The equations defining $Y = X \times_{V_1} W_1$ 
    in $\mathbb C^n \times \mathbb C \times \mathbb P^{t-1}$
    are again
    \[ M(\underline x, \varepsilon) \cdot \vec s = 0 \] 
    with $\vec s = \left( s_1, \dots, s_t \right)^T$ the vector whose entries
    are the homogeneous coordinates of $\mathbb P^{t-1}$. 
    For the special fiber $Y_0 = X_0 \times_{V_1} W_1$ one always obtains
    the same result as in the case of applying the Tjurina modification to 
    the singularity
    alone by Construction.
    \label{cst:tjurinaModificationInFamily}
\end{construction}

\begin{example} 
        \label{exp:3by4MatrixInFamily}
   Consider the deformation with a parameter $\varepsilon$ given 
   by the matrix
$$
    M(\underline x, \varepsilon) = 
     \begin{pmatrix}
        x & y-v & y+z+2\varepsilon \cr
        y & z-v & x+u+2\varepsilon \cr
         z & 3\varepsilon & x-u \cr
        3\varepsilon & u & v
    \end{pmatrix}.
$$
Let ${X} \subset \mathbb C^5 \times \mathbb C$ be the total space 
of the deformation
\[
    \begin{xy}
        \xymatrix{
            X_0 \ar[r] \ar[d] & X \ar[d]^\varepsilon \\
            \{0\} \ar[r] & \mathbb C
        }
    \end{xy}
\]
The Tjurina modification in family 
in $\mathbb C^5 \times \mathbb C \times \mathbb P^2$
is now described by the equations
\[ M(\underline x, \varepsilon) \cdot 
    \begin{pmatrix}
        s_1 \\ s_2 \\ s_3
    \end{pmatrix} 
    = 0
\]
As a direct computation shows, all fibers of this family except the 
one over $\varepsilon=0$ are smooth. This has an important consequence
in the setting of Tjurina modification: 
\end{example}

\begin{proposition}
    If in the setting of Construction
    \ref{cst:tjurinaModificationInFamily} above the 
    deformation of $X_0$ over $\mathbb C$ is a smoothing, 
    the restriction of $\pi$ to a smooth fiber 
    \[
        \pi_\varepsilon : Y_\varepsilon \to X_\varepsilon 
    \]
    in diagram (\ref{eqn:DiagramTjurinaModInFamily}) 
    is an isomorphism.
    \label{prp:IsomorphismSmoothFibersTjurinaModification}
\end{proposition}

\begin{proof}
    For fixed $\varepsilon$ the rational map 
    \[ P \circ M( -,\varepsilon) : X_\varepsilon \dashrightarrow \mathbb P^{t-1} \]
    is not well defined in the vanishing locus of $(t-1)$-minors of $M$.
    By Corollary \ref{cor:MaxMinorsOrigin}, this is contained in the singular 
    locus of $X_\varepsilon$, which is empty for smooth fibers.
    Hence $P \circ M_\varepsilon$ is regular.
\end{proof}

Unfortunately it is not at all clear that 
the families $Y \overset{\varepsilon \circ \pi}{\longrightarrow} \mathbb C$
obtained by Tjurina modifications in family as in 
diagram (\ref{eqn:DiagramTjurinaModInFamily}) are flat. 
Whether or not this is the case, will in general depend on the 
deformation in question, the dimension of the singularity and the 
Cohen-Macaulay-type. 
Consider e.g. a space curve $X_0 \subset \mathbb C^3$ of Cohen-Macaulay-type
$3$ as the special fiber in a smoothing by a parameter 
$\varepsilon$, then the fiber $Y_0$ of the Tjurina transform over 
$\varepsilon = 0$ contains a $\mathbb P^2$, while the other fibers stay 
$1$-dimensional. This clearly contradicts flatness. 

For simple ICMC2 singularities of dimension $\dim(X_0,0)>0$ this does
not pose a problem\footnote{We deliberately exclude the simple fat points
here as their behaviour obviously differs from the higher dimensions, because
any ${\mathbb P}^k$ in the $Y_0$ would violate flatness.}:  All families 
in the classification of Neumer and Fr\"uhbis-Kr\"uger \cite{FN} have 
Cohen-Macaulay type $t = 2$. This will turn out to be sufficient to 
assure flatness in all of our cases of interest.

\begin{proposition}
    Let $(X_0,0)\subset (\mathbb C^{n+2},0)$ 
    be an ICMC2 singularity of dimension $n > 0$ and Cohen-Macaulay type 
    $t \leq n+1$.
    The Tjurina modification in family for a deformation 
    $X_0\hookrightarrow X \overset{\varepsilon}{\longrightarrow} \mathbb C$ 
     \[
        \begin{xy}
            \xymatrix{
                Y_0 \ar[d]_{\pi_0} \ar[r] & Y \ar[d]^\pi \\
                X_0 \ar[r] \ar[d] & X \ar[d]^\varepsilon \\
                \{0\} \ar[r] & \mathbb C
            }
        \end{xy}
    \]
    is flat over $\mathbb C$.
    \label{prp:FlatnessTjurinaModificationSimpleSingularities}
\end{proposition}

\begin{proof}
    The Grassmannian in question is a $\mathbb P^{t-1}$.
    As usual let $(s_1:\dots:s_t)$ be its projective coordinates 
    and $M(\underline x, \varepsilon)$ 
    the matrix describing the family $X$.
    The variety $Y_0 \subset \mathbb C^{n+2} \times \mathbb P^{t-1}$
    is given by the $t+1$ equations 
    \[
        M_0(\underline x) \cdot 
        \begin{pmatrix}
            s_1 \\ \vdots \\ s_t
        \end{pmatrix}
        = 0
    \]
    Now 
    \[
        \dim Y_0 = \max \left\{ \dim X_0, \dim \mathbb P^{t-1} \right\} 
        = \max \left\{ n, t-1 \right\} = \dim X_0,
    \]
    because $\pi_0 : Y_0 \to X_0$ is an isomorphism on 
    $Y_0 \setminus \pi_0^{-1}(\{0\})$ 
    and the exceptional set $\pi_0^{-1}(\{0\})$ is a $\mathbb P^{t-1}$.
    Since $X_0$ had codimension $2$ in $\mathbb C^{n+2}$ we find 
    $Y_0$ to have codimension $t+1$ in 
    $\mathbb C^{n+2} \times \mathbb P^{t-1}$.
    But locally in all charts, there are exactly $t+1$ equations describing
    $Y_0$.
    This means $Y_0$ is a locally complete intersection, 
    so the induced deformation by $M(\underline x,\varepsilon)$ in 
    the Tjurina modification in family is flat.
\end{proof}

\begin{remark}
    The above result was independently formulated by Jesse Kass 
    for simple space curve singularities in his up-to-now unpublished
    work on Coxeter-Dynkin diagrams of space curve singularities. 
\end{remark}

An alternative way to check flatness of a family is checking the relation
lifting property for the relations of generators of the defining ideal
(cf. \cite{Art}).

Let $X_0 \subset \mathbb C^n$ be an ICMC2 singularity of Cohen-Macaulay type 
$t$ and $M(\underline x, \varepsilon)$ be a matrix defining a deformation of 
$X_0$ over $(\mathbb C,0)$. 
The ideal $J \subset \mathbb C\{\underline x,\varepsilon\}[s_1,\dots,s_t]$
defining the Tjurina transform 
$Y \subset \mathbb C^n \times \mathbb C \times \mathbb P^{t-1}$ 
is generated by the $t$ equations $H_i(\underline x, \varepsilon, \underline s) = 0$
originating from the lines of the system
\[ M(\underline x, \varepsilon ) \cdot 
    \begin{pmatrix}
        s_1 \\ \vdots \\ s_t
    \end{pmatrix}
    = 0.
\]
The relation lifting property for flatness requires that any relation 
\[ \sum_j r_j \cdot h_j = 0\]
in $\mathbb C\{\underline x\}[\underline s]$ among the 
$h_j = H_j(\varepsilon = 0)$ can be lifted to a relation 
$\sum_j R_j \cdot H_j = 0$ in 
$\mathbb C\{\underline x, \varepsilon\}[\underline s]$
with $r_j = R_j(\varepsilon = 0)$.
Now there is one relation among the $H_j$ which comes naturally 
with a lifting:

Because the matrix $M$ describes the syzygies of the generators of 
\[
    I = \langle F_1(\underline x, \varepsilon), \dots, F_{t+1}
    (\underline x, \varepsilon) \rangle,
\]
i.e. the ideal defining 
$X \subset \mathbb C^n \times \mathbb C$, 
we can write 
\[
    0 = 
    \left( F_1,\dots, F_{t+1} \right) \cdot M \cdot 
    \begin{pmatrix}
        s_1 \\ \vdots \\ s_t
    \end{pmatrix}
    = 
    \left( F_1,\dots, F_{t+1} \right) \cdot 
    \begin{pmatrix}
        H_1 \\ \vdots \\ H_t
    \end{pmatrix}.
\]
We call this the ``relation by the maximal minors''. 
This leads to the following criterion for flatness of the Tjurina modification.

\begin{lemma}
    If in the above setting the relations among the generators $h_i$ of 
    the ideal defining the Tjurina modification $Y_0$ of an ICMC2 singularity
    $X_0$ are generated by the Koszul relations and the relation by the 
    maximal minors, then any Tjurina modification in family of a deformation
    of $X_0$ is again flat.
    \label{lem:FlatnessRelations}
\end{lemma}

We finish the discussion of flatness by looking at 
one last example which does not satisfy the condition in the preceding 
criteria:

\begin{example}
The family of ICMC2 fat points defined by
$$
\begin{pmatrix}
          0 & x \cr
          x & y \cr
          y & \varepsilon 
\end{pmatrix}
$$
does not give rise to a flat family by Tjurina modification:\\
Tjurina modification provides 
$$\hat X_\varepsilon = V(\underbrace{s_2x}_{=f_1},
                         \underbrace{s_1x+s_2y}_{=f_2},
                         \underbrace{s_1y}_{=f_3}+s_2\varepsilon)$$
For $\varepsilon=0$ we have the additional relation
$$s_1^2 f_1 - s_1s_2 f_2 + s_2^2 f_3 =0$$
among the generators of the ideal of $Y_0$. The relation cannot be
lifted to a relation of the whole family.\\
Considering this example from a geometric perspective, all fibers 
except the fiber at $\varepsilon=0$ are zero-dimensional, but the special
fiber additionally contains the $\mathbb P^1$ introduced by the Tjurina
modification. As before such a jump in dimension clearly contradicts flatness.
\end{example}

To end this section, we want to study the effects of a Tjurina modification 
to versal families of ICMC2 singularities of Cohen-Macaulay type $t=2$.
These observations originate from direct computations, but will be  
useful for explicit examples:

\begin{remark}
Let $(X_0,0) \subset (\mathbb C^n,0)$ be an ICMC2 singularity at the origin 
with Cohen-Macaulay type $t=2$ and $\dim X_0 >0$. Let
$M$ be the corresponding presentation matrix. 
Expanding the matrix entries up to degree $r$ and taking equivalence classes
modulo $\langle x_1,\ldots,x_n\rangle^{r+1}$, we can represent each such 
class by a matrix with polynomial entries of degree at most $r$. We shall
refer to this representative as the $r$-jet of the presentation matrix, 
$j_r M$.
More precisely, we need to prepare the subsequent discussion of the 
relationship of the deformations of $X_0$ and $Y_0$ and thus determine 
a very coarse classification of occurring $1$-jets\footnote{This list
is, of course, loosely related to the lists of $1$-jets in \cite{FN}, 
but it contains significantly fewer classes, because here we are only 
hunting for a criterion for isolatedness of the singularities of the 
Tjurina transform.}.

As we are considering a germ around the origin, all entries of 
$j_1 M$ are homogeneous linear polynomials. We know that row and column
operations on $M$ leave the germ $(X_0,0)$ unchanged, and we can safely 
pass to sufficiently general ${\mathbb C}$-linear combinations of the 
two original columns. The second column of $j_1 M$
thus holds up to $3$ ${\mathbb C}$-linearly independent linear forms.
By suitable row operations on $M$, we can then cancel linearly dependent
entries of this column of $j_1 M$ and achieve that the zero entries 
are positioned below the non-zero entries. (Note that the sufficiently 
general linear combination of the columns now ensures that a row with
a zero in the second entry also holds a zero in the first entry.)  
By an analytic change of coordinates, we can now choose the non-zero
entries of the second column as new coordinates, starting with $x_1$,
and obtain the following four cases: 
$$\begin{pmatrix} * & x_1 \cr * & x_2 \cr * & x_3\end{pmatrix},\;\;\;
\begin{pmatrix} * & x_1 \cr * & x_2 \cr 0 & 0\end{pmatrix},\;\;\;
\begin{pmatrix} * & x_1 \cr 0 & 0 \cr 0 & 0\end{pmatrix},\;\;\;
\begin{pmatrix} 0 & 0 \cr 0& 0\cr 0 & 0\end{pmatrix}.$$
Here $*$ denotes an arbitrary entry. 
\end{remark}

\begin{lemma} 
    \label{lem:ISinTjurinaTransfToJetOfMatrix}
Let $(X_0,0) \subset (\mathbb C^5,0)$ be an ICMC2 threefold 
singularity of Cohen-Macaulay type $t=2$.
Then the Tjurina transform 
$Y_0$ has at most isolated singularities, iff $X_0$ is 
contact equivalent to a ICMC2 with presentation matrix 
$$\begin{pmatrix} a & x_1 \cr b & x_2 \cr c & x_3 \end{pmatrix}, $$
where $a,b,c \subset \langle x_1,\dots,x_5\rangle$.
\end{lemma}

\begin{proof}
Tjurina 
modification is an isomorphism outside the singular locus,  which implies 
that the singular locus of $Y_0$ is contained in 
$E = \pi_0^{-1}(\{0\}) \cong \mathbb P^1$.
Because $E$ is irreducible, the singular locus of $Y_0$ is either a 
finite number of points or the whole $\mathbb P^1$. \\
If the matrix has the desired structure, we focus on one of the two 
standard affine charts of the exceptional ${\mathbb P}^1$ to show that there
can be at most isolated singularities. 
As usual let $(s_1:s_2)$ be the homogeneous coordinates of the 
exceptional curve $\mathbb P^1$. 
The equations for $Y_0 \subset \mathbb C^5 \times \mathbb P^1$ 
in the chart $s_1 \neq 0$ are given by the ideal
$I= \langle a+s_2x_1, b+s_2x_2, c+s_2x_3 \rangle$.
The jacobian of this complete intersection reads
\begin{equation}
\begin{pmatrix} \frac{\partial a}{\partial x_1} + s_2 &
                  \frac{\partial a}{\partial x_2} &
                  \frac{\partial a}{\partial x_3} & 
                  \frac{\partial a}{\partial x_4} & 
                  \frac{\partial a}{\partial x_5} & x_1\cr
                  \frac{\partial b}{\partial x_1} &
                  \frac{\partial b}{\partial x_2} + s_2 &
                  \frac{\partial b}{\partial x_3} &
                  \frac{\partial a}{\partial x_4} & 
                  \frac{\partial b}{\partial x_5} & x_2 \cr
                  \frac{\partial c}{\partial x_1} &
                  \frac{\partial c}{\partial x_2} &
                  \frac{\partial c}{\partial x_3} + s_2 & 
                  \frac{\partial a}{\partial x_4} & 
                  \frac{\partial c}{\partial x_5} & x_3\cr
\end{pmatrix}.
    \label{eqn:JacobianG}
\end{equation}
One of its 3-minors (first 3 columns) and hence of the ideal of the 
singular locus contains an element of the form $s_2^3+\phi$ where the 
$s_2$-degree of the remaining part $\phi$ is at most $2$. This excludes 
the case of the singular locus being the whole exceptional curve $E$.

If, on the other hand, the matrix is not of the desired form, at least one
row and hence at least one generator of $I$ is contained in 
$\langle x_1,\dots,x_5 \rangle^2$, whence at least one row of the jacobian
matrix -- and thus the ideal of its $3$-minors -- is contained in 
$\langle x_1,\dots,x_5 \rangle$. Hence the singular locus would be 
$1$-dimensional in this case.
\end{proof}

In the case of the preceding lemma with only isolated singularities 
in the Tjurina transform, we 
now compare the infinitessimal deformations of the ICMC2 
singularity $(X_0,0)$ 
downstairs with those of the local complete intersection scheme $(Y_0,\PP^1\times\{0\})$ 
upstairs, where $\PP^1 \times \{0\}$ is the exceptional locus of 
the Tjurina transform.
For the affine germ $(X_0,0)$ the first order deformations 
are encoded in the $\CC\{\underline x\}$-module 
$T^1_{X_0}$. 
The space of embedded first order deformations for the 
Tjurina transform $\iota : (Y_0,\PP^1\times \{0\}) \hookrightarrow \PP^1\times \CC^5$ 
can be described as follows. 
Let $\mathcal I$ be the ideal sheaf defining 
$(Y_0,\PP^1 \times \{0\})$ in $(\PP^1\times \CC^5, \PP^1\times \{0\})$. 
We take global sections of the normal bundle 
\[
    N_{Y_0} = H^0(Y_0,\sheafhom_\mathcal{O}(\mathcal{I},\mathcal O_{Y_0}))
\]
and divide by those deformations coming from global sections 
of the tangent bundle $H^0(Y_0,\iota^* T_{\PP^1\times \CC^5})$. 
The resulting quotient will be denoted by
\begin{equation}
    N' := N_{Y_0}/H^0(Y_0,\iota^*T_{\PP^1\times \CC^5}).
    \label{eqn:FirstOrderDefY0}
\end{equation}
Note that the global section functor takes coherent sheaves to 
finitely generated $\CC\{\underline x\}$-modules.
In fact $N'$ is naturally a $\CC\{\underline x\}$-module with 
support in the point $0$ and hence a finite dimensional 
vector space over $\CC$. To see this observe that outside the 
singular locus $0\in X_0$ (and outside $\PP^1\times \{0\} \subset Y_0$ respectively), 
the space $Y_0$ is described as a graph over $X_0$ and we therefore have 
a natural splitting of the normal bundle 
\[
    N_{Y_0} = N_{X_0} \oplus T_{\PP^1}|_{Y_0}.
\]
Because the tangent bundle of $\PP^1$ is globally generated, 
the second summand is killed when forming the quotient $N'$. 
But the first summand cancels on the smooth locus anyway. 

It is clear from the construction that every deformation of $(X_0,0)$ 
induces a deformation of $(Y_0,\PP^1 \times \{0\})$. Let $(X_0,0)$ be given  
by the matrix 
\[
    M = 
    \begin{pmatrix}
        a & x_1 \\
        b & x_2 \\
        c & x_3
    \end{pmatrix}
    \in \Mat(3,2; \mathbb C\{x_1,\dots,x_5\}).
\]
and let 
\[
    H_1 = s_1 \cdot a + s_2 \cdot x_1,\quad 
    H_2 = s_1 \cdot b + s_2 \cdot x_2, \quad
    H_3 = s_1 \cdot c + s_2 \cdot x_3  
    \in \mathbb C\{\underline x\}[s_1,s_2]
\]
be the three equations defining the Tjurina transform $Y_0$ in 
$\PP^1 \times \mathbb C^5$, which are homogeneous in $\underline{s}$.
On the level of equations there is a map 
\begin{equation}
    \begin{xy}
        \xymatrix{
            \left\{
                \textnormal{Perturbations of } M 
            \right\} \ar[r]^\Lambda \ar@{-}[d]_{1:1} & 
            \left\{
                \textnormal{Perturbations of } \underline H
            \right\} \ar@{-}[d]^{1:1} \\
            \Mat(3,2;\mathbb C\{\underline x\})\ar[r]^\Lambda_{E_{i,j}^{(2,3)} \mapsto e_is_j} & 
            \left( (\mathbb C\{\underline x\}[s_1,s_2])^3 \right)_{(1)}, 
        }
    \end{xy}
    \label{eqn:DefinitionLambda}
\end{equation}
where the $e_i$ denote the generators of the free module on the right hand side
and $E_{i,j}^{(r,s)}$ denote the 
$r \times s$ matrices possessing only one non-zero entry of value $1$ at 
position $i,j$. The lower index $(1)$ signifies that we only consider the 
homogeneous part of degree $1$ in $\underline{s}$. 

\begin{lemma}
    The map $\Lambda$ induces an isomorphism of 
    first order deformations of $(X_0,0)$ and $(Y_0,E)$, 
    i.e. an isomorphism of 
    $\CC\{\underline x\}$-modules 
    \[
        \Lambda: T^1_{X_0,0} \overset{\cong}{\longrightarrow} N'.
    \]
    \label{lem:LambdaIsomorphismus}
\end{lemma}

\begin{proof}
    We have already obtained the isomorphism $\Lambda$ between 
    $\Mat(3,2;\mathbb C\{\underline x\})$ and 
    $\left( \mathbb C\{\underline x\}[s_1,s_2] \right)^3_{(1)}$.
    From the description of the $T^1_{X_0,0}$ in Lemma \ref{lem:structT1} and the
    definition of $N'$ we know the relations on both sides.
    It hence remains to prove that 
    the modules $J_M + \operatorname{Im}(g)$ from Lemma \ref{lem:structT1} and
    $(J_H + I_H)_{(1)}$ are isomorphic.
    Here 
    $I_H= \langle H_1,H_2,H_3 \rangle \mathbb C\{\underline x\}^3$ and
    $J_H$ is generated by the columns of the Jacobian matrix of the $H_i$ 
    defining $Y_0$.\\
    By construction of $\underline{H}$, we see immediately
    $$ \Lambda(\frac{\partial M}{\partial x_i}) = 
       \frac{\partial \underline{H}}{\partial x_i},$$
    $$ \Lambda(M \cdot E^{(2,2)}_{i,j}) = 
       s_i \frac{\partial \underline{H}}{\partial s_j}$$
    and
    $$ \Lambda(E_{i,j}^{(3,3)} \cdot M) = H_j e_i.$$
    This provides a $1:1$ correspondence of the generators of these two modules
     and hence proves the claim about the cokernels:
    \begin{equation*}
    \xymatrix{
    0 \ar[r] & J_M + \operatorname{Im}(g) \ar[r] 
                                          \ar[d]^{\cong}_{\Lambda} 
             & \Mat(3,2;\mathbb C\{\underline x\}) \ar[r] 
                                                   \ar[d]^{\cong}_{\Lambda}
             & T^1_{X_0,0} \ar[r] \ar[d]
             & 0\\
    0 \ar[r] & (J_H + I_H)_{(1)} \ar[r]
             & \left( \mathbb C\{\underline x\}[s_1,s_2] \right)^3_{(1)} \ar[r]
             & N' \ar[r]
             &0
    }
    \end{equation*}
\end{proof}

There is a splitting of the module $N'$ coming from the 
local-to-global spectral 
sequence of the exact sequence of sheaves
\begin{equation}
    \xymatrix{
        0 \ar[r] &
        T_{Y_0} \ar[r] &
        \iota^* T_{\PP^1\times \CC^5} \ar[r] & 
        N_{Y_0} \ar[r] &
        T^1_{Y_0} \ar[r] &
        0,
    }
    \label{eqn:ExactSequenceT1Y0}
\end{equation}
which can be explicitly described as follows. 

We first split the exact sequence (\ref{eqn:ExactSequenceT1Y0}) into 
short exact sequences
\begin{equation}
    \xymatrix{
        0 \ar[r] &
        T_{Y_0} \ar[r] &
        \iota^* T_{\PP^1\times \CC^5} \ar[r] & 
        \mathcal K \ar[r] &
        0 & & \\
        & & 0 \ar[r] & 
        \mathcal K \ar[r] &
        N_{Y_0} \ar[r] &
        T^1_{Y_0} \ar[r] &
        0\\
    }.
    \label{eqn:ExactSequenceT1Y0Splitting}
\end{equation}
The long exact sequences in cohomology both have to finish after 
the degree one terms, because the underlying scheme is covered by 
two affine charts. 

Let again $\mathcal I$ be the ideal sheaf of $(Y_0,\PP^1)$. If we 
tensor the short exact sequence 
\[
    \xymatrix{
        0 \ar[r] &
        \mathcal I \ar[r] &
        \mathcal O_{\PP^1\times \CC^5} \ar[r] & 
        \mathcal O_{Y_0} \ar[r] &
        0
    }
\]
with the locally free sheaf 
$T_{\PP^1 \times \CC^5}$ and take the long exact sequence in 
cohomology, we see that 
\[
    H^1(Y_0,\iota^* T_{\PP^1\times \CC^5}) = 0.
\]
Looking at the first long exact sequence in cohomology of 
(\ref{eqn:ExactSequenceT1Y0Splitting}), we deduce that 
\begin{equation}
    \coker\left( H^0(Y_0, \iota^* T_{\PP^1 \times \CC^5}) \to H^0(Y_0,\mathcal K) \right) 
    \cong H^1(Y_0,T_{Y_0)}) 
    \label{eqn:GloballyNontrivialDeformationsH1}
\end{equation}
and 
\begin{equation}
    H^1(Y_0,\mathcal K) = 0.
    \label{eqn:VanishingH1CokernelTangentbundle}
\end{equation}
Combining these results with the second long exact sequence 
of (\ref{eqn:ExactSequenceT1Y0Splitting}) and recalling that  
$N'=N_{Y_0}/H^0(Y_0, \iota^* T_{\PP^1 \times \CC^5})$, we obtain a 
short exact sequence 
\[
    \xymatrix{
        0 \ar[r] &
        H^1(Y_0,T_{Y_0}) \ar[r] &
        N' \ar[r] &
        H^0(Y_0,T^1_{Y_0}) \ar[r] &
        0
    }, 
\]
the middle term of which is a finite dimensional vector space over 
$\CC$. Any choice of a splitting gives us 
\begin{equation}
    N' = 
    H^1(Y_0,T_{Y_0}) \oplus H^0(Y_0,T^1_{Y_0}).
    \label{eqn:DecompositionNprime}
\end{equation}
The sheaf underlying the right hand side summand is supported 
only in the singular points and hence affine. Thus if we let 
$\Sigma(Y_0)$ be the set of singular points of $Y_0$ we can 
rewrite (\ref{eqn:DecompositionNprime}) as 
\begin{equation}
    N' = 
    H^1(Y_0,T_{Y_0}) \oplus \bigoplus_{p\in \Sigma(Y_0)} T^1_{Y_0,p}
    \label{eqn:DecompositionNprime2}
\end{equation}
In particular for any $q \in \Sigma(Y_0)$ we get a surjective map 
from $T^1_{X_0,0}$ onto $T^1_{Y_0,q}$ by the composition
\[
    T^1_{X_0,0} \cong N' \cong 
    H^1(Y_0,T_{Y_0}) \oplus \bigoplus_{p\in \Sigma(Y_0)} T^1_{Y_0,p} 
    \longrightarrow 
    T^1_{Y_0,q},
\]
where the last map is the projection to the summand for $q$.
This proves the following corollary.

\begin{corollary}
    \label{cor:VersalityDefTjurinaTransf}
    Let $(X_0,0)\subset (\mathbb C^5,0)$ be an ICMC2 threefold singularity of 
    Cohen-Macaulay type $t=2$ such that the Tjurina transform 
    $Y_0$ has at most isolated singularities.
    Furthermore let
    $X_0 \hookrightarrow X \longrightarrow {\mathbb C}^{\tau}$ be a 
    semi-universal deformation of $X_0$. Then the induced family 
    $Y_0 \hookrightarrow Y \longrightarrow {\mathbb C}^{\tau}$ is again 
    versal for each of the arising singularities.
\end{corollary}

Note that the induced local deformations for the isolated singularities
of $Y_0$ do not need to be semi-universal, i.e. $\tau$ might not 
be minimal.

\begin{remark}
    As can be expected given the results of this section, for all $(X_0,0)$ 
    in the table of simple ICMC2 singularities of dimension
    $3$, the Tjurina modification $Y_0$ has at most simple ICIS as can be
    read off from table \ref{tab:ResultsThreefolds}.
    \label{rem:SimpleDownSimpleUp}
\end{remark}

%% file: chapter4.tex
\section{Vanishing cycles}

From now on we'll often be concerned with the homology groups of a 
given topological space. By this we mean Simplicial homology 
with integer coefficients and we'll just write $H_\bullet(-)$
for $H_\bullet(-,\mathbb Z)$ for short.
To fix notation, we briefly recall the definition of Milnor fiber and 
vanishing cycles for isolated singularities (see e.g. \cite{Mil68} for
a reference on these topics).

Let $(X_0,0) \subset (\mathbb C^n,0)$ be an isolated singularity
at the origin and $X_0$ a representative thereof. Then there is a real $\eta_0 > 0$ 
such that the intersection of $X_0$ with the sphere $S_\eta$ 
of radius $\eta$ is transversal for all $\eta_0 \geq \eta > 0$.
For any $\eta >0$ chosen in this way, we will refer to a closed ball $B_\eta$ 
of radius $\eta$ around $0 \in \mathbb C^n$ as a Milnor ball for the 
singularity $X_0$. Furthermore we denote by $\overline{X_0}$ the topological space 
\begin{align*}
    \overline X_0 := X_0 \cap B_\eta.
\end{align*}
We explicitly cite the following well-known theorem, which ensures that 
these definitions are independent of the chosen (sufficiently small) $\eta$.

\begin{theorem}[conical structure, \cite{Mil68}]
    For an isolated singularity $(X_0,0) \subset (\mathbb C^n,0)$ and 
    a Milnor ball $B_\eta$ for $X_0$ the pair of spaces 
    $(B_\eta, \overline X_0)$ is homeomorphic to the pair
    $(C(S_\eta), C(\partial \overline X_0))$, 
    where $C(L)$ denotes the cone over $L\subset S_\eta$,
    i.e. the set of real line segments to the origin.
    This can be chosen to be a diffeomorphism on the open set 
    $(B_\eta \setminus \{ 0 \}, \overline X_0 \setminus \{0\})$.
    \label{thm:conicalStructure}
\end{theorem}

Consequently $\overline X_0$ and $\partial\overline X_0$ are 
well defined topological spaces up to homeomorphism
for a germ $(X_0,0)$ of an isolated singularity, i.e. do not depend
on the representative $X_0$.
Now, consider a deformation of an isolated singularity
$(X_0,0) \subset \mathbb (\mathbb C^n,0)$ by some parameter $\varepsilon$, 
i.e. a flat family
    \begin{align*}
        \begin{xy}
            \xymatrix{ 
                X_0 \ar[d] \ar[r]^\iota & X \ar[d]^\varepsilon \\
                \{ 0 \} \ar[r] & \mathbb C
            }
        \end{xy}
    \end{align*}
where $X_0 \subset \mathbb C^n$ and $X\subset \mathbb C^n \times \mathbb C$ 
are representatives of the respective germs.
Having chosen a Milnor ball $B_\eta$ for $X_0$ there exists
an open neighborhood $0 \in D \subset \mathbb C$ in the 
deformation base $\mathbb C$, such that for all $\varepsilon \in D$
the intersection 
    \[
        \partial \overline X_\varepsilon = 
         X_\varepsilon \cap \partial B_\varepsilon
    \]
of the fiber $X_\varepsilon$ with the Milnor ball in the 
fiber over $\varepsilon$ is transversal. The cylinder $B_\eta \times D$ 
is called a Milnor tube for the deformation of $X_0$. \\

Theorem \ref{thm:conicalStructure} ensures that all the homology
groups of $\overline X_0$ vanish except in degree $0$. 
But in any deformed fiber $\overline X_\varepsilon$ there may 
exist nontrivial cycles. 
If a fiber $\overline X_\varepsilon$ is smooth, i.e. a smooth complex 
manifold with boundary, it is called a Milnor fiber of the singularity $X_0$.
Any nontrivial cycles in the homology of $\overline X_\varepsilon$
are called vanishing cycles of the singularity $X_0$.
It is well known what these vanishing cycles look like for 
any ICIS of complex dimension $n$: they form a 
bouquet of spheres of real dimension $n$, see \cite{Hamm}

For the ICMC2 singularities of dimension 3 which we are considering in this 
article, known results on the vanishing cycles of the Milnor fibers 
are scarce.

\begin{remark}
    It is a priori not clear and in general wrong 
    to expect exactly one Milnor fiber for a given isolated singularity 
    $X_0$ (up to diffeomorphism). First of all 
    there may not exist any deformation with smooth fibers at all, like
    for the rigid isolated 4-fold singularity appearing in the classification
    of Fr\"uhbis-Kr\"uger and Neumer.\\
    On the other hand, given two different smoothings 
    $\pi : X \to \mathbb C$ and $\pi' : X' \to \mathbb C$ two smooth
    fibers $X_\varepsilon$ and $X'_\varepsilon$ are not nescessarily
    diffeomorphic as the famous example of Pinkham \cite{Pin} shows. 
    They are, however, diffeomorphic if they belong to the same connected 
    component of the deformation base. This is an immediate corollary 
    of the Ehresmann fibration theorem. 
\end{remark}

If $X_0$ is a smoothable ICMC2 singularity, the set of points 
$\varepsilon \in \mathbb C^\tau$ with smooth fibers is open and 
connected since its complement (the discriminant) has real codimension 
at least $2$. Therefore in all our cases of interest in this article
the singularities have a unique Milnor fiber up to diffeomorphism.

We need one more preliminary result which will be applied to determine
the topology of the Tjurina modification $Y_0$.

\begin{proposition}
    Let $(X_0,0) \subset (\mathbb C^n, 0)$ be an isolated singularity
    and $ \pi_0 : Y_0 \to X_0$ a morphism defined on suitably small 
    representatives such that the restriction
    \[
        \pi_0 : Y_0\setminus \pi^{-1}(\{0\}) \to X_0 \setminus \{0\}
    \]
    is an isomorphism and the exceptional set 
    $E = \pi_0^{-1}(\{0\})$ is closed and projective. 
    Then $E$ is a deformation retract of $Y_0$.
    \label{prp:ExceptionalSetENR}
\end{proposition}

\begin{proof}
    The variety $E$ is closed and projective, hence compact. 
    It follows from 
    \cite{Loj}, 
    that $E$ is a Euclidean Neighborhood Retract of an open 
    neighborhood $U$ of $E$ in $Y_0$. 
    But outside $E$ the map $\pi_0$ is an isomorphism, so
    $\pi_0(U) \subset X_0$ is open.
    With the theorem about the conical structure 
    (\ref{thm:conicalStructure}) 
    we can now shrink 
    $\overline Y_0 \setminus E = \overline X_0 \setminus\{0\}$
    to something homotopic to $\overline Y_0$ inside the open set 
    $\pi_0(U)$ 
    and subsequently to $E$.
\end{proof}

Using Tjurina modification in family, we are now ready to explain 
the observations of \cite{DP2} in the case of a simple ICMC2 threefold 
$(X_0,0) \subset (\mathbb C^5,0)$. 
Applying the Tjurina modification we get a transform $Y_0$ with only
A-D-E singularities according to Remark \ref{rem:SimpleDownSimpleUp}. 
Since we nescessarily have Cohen-Macaulay type $t=2$, 
the homotopy type of $Y_0$ is given by the exceptional set 
$\mathbb P^1 = E \subset Y_0$ as a consequence of Proposition 
\ref{prp:ExceptionalSetENR}. So we always find
\begin{equation}
    b_0(Y_0) = 1,\quad b_1(Y_0) = 0,\quad b_2(Y_0) = 1,\quad b_3(Y_0) = 0.
    \label{eqn:BettiNumbersTjurinaTransfSimpleSing}
\end{equation}
Now Proposition \ref{prp:FlatnessTjurinaModificationSimpleSingularities} 
assures the Tjurina modification to be well behaved within families. 
Hence we can choose any smoothing 
$ X_0 \hookrightarrow X \overset{\varepsilon}{\longrightarrow} \mathbb C$
and carefully observe the interplay of cycles present in $Y_0$ 
with upcoming vanishing cycles of the ICIS when passing 
from $Y_0$ to a deformed fiber $Y_\varepsilon$ in the induced deformation.
This is covered in Theorem \ref{thm:TjTransfInterplay}. 
Finally we can use the identification $Y_\varepsilon \cong X_\varepsilon$ 
from Proposition \ref{prp:IsomorphismSmoothFibersTjurinaModification}
to obtain the desired vanishing topology.\\

We slightly weaken the assumptions 
and also allow $(X_0,0)\subset (\mathbb C^5,0)$ to be a 
non-simple ICMC2 threefold singularity.
However we still require the Cohen-Macaulay type to be $t=2$ and 
the jet type as in Lemma \ref{lem:ISinTjurinaTransfToJetOfMatrix}, 
i.e. only ICIS in $Y_0$; the more general case allowing non-isolated 
singularities in the Tjurina transform will be studied in \cite{Z1} 
in more detail. 
The results are gathered in the Tables \ref{tab:ResultsThreefolds}
and \ref{tab:boundingThreefolds} in the next section.

\begin{theorem} \label{thm:TjTransfInterplay}
    In the above setting 
    consider $X_0$ as the special fiber in a smoothing
    $X_0 \hookrightarrow X \overset{\varepsilon}{\longrightarrow} \mathbb C$
    together with the Tjurina modification in family
    $Y_0 \hookrightarrow Y 
         \overset{\varepsilon\circ \pi}{\longrightarrow} \mathbb C$
    as in diagram \ref{eqn:DiagramTjurinaModInFamily}.
    We denote the Milnor tube arising from $B$ by 
    \[ T = B \times D \subset \mathbb C^5 \times \mathbb C \]
    and the one originating from $\hat B$ by
    \[
        \hat T = \pi^{-1}( T )\subset 
                 \mathbb C^5 \times \mathbb C \times \mathbb P^1,
    \]
    which also allows us to refer to the Milnor fiber
    $\overline X_\varepsilon$ and the fiber 
    $\overline Y_\varepsilon = \hat T \cap Y_\varepsilon$ sitting over it.
    The Betti numbers of a smooth fiber $\overline X_\varepsilon$ are given by
    \[
        b_0 = 1, \quad b_1 = 0, \quad b_2 = 1, \quad b_3 = r 
    \]
    where $r \in \mathbb N$ is the sum of the Milnor numbers of the 
    ICIS of $Y_0$.
\end{theorem}

\begin{proof}
    Throughout the proof many steps require shrinking the open set 
    of admissible deformation parameters $\varepsilon$ in the deformation
    base. However, those steps are finitely many and no harm is done, since 
    we only consider representatives of germs. For the reader's convenience, 
    we will suppress mentioning this obvious technical detail each time 
    it occurs.

    Let $n$ be the number of singularities of $\overline Y_0$. 
    Fix local analytic embeddings of these ICIS 
    to some affine space and let
    $T = \bigcup_{i=1}^n T_i$ be a collection of Milnor tubes 
    around the singularities of $\overline Y_0$ for the induced 
    deformation 
    $\overline Y_0 \hookrightarrow \overline Y 
    \overset{\varepsilon\circ \pi}{\longrightarrow} \mathbb C$.
    For each $i=1,\dots,n$, let $S_i$ 
    be the Milnor tubes sweeped 
    out by Milnor balls of half of the radius of $T_i$.
    Decompose the total space $\overline Y$ in the two open sets 
    \[ U := \overline Y \setminus \left( \bigcup_{i=1}^n S_i \cap Y \right) \quad
        \textnormal{ and } \quad
        V \textnormal{ the interior of }\left( \bigcup_{i=1}^n T_i \cap Y \right)
    \]
    Let $\overline U, \overline V$ and $\overline W=\overline{U \cap V}$ be 
    the closures of the respective open sets in the Euclidean topology. Each 
    of them is compact, $\overline U$ and $\overline W$ even compact manifolds
    with boundary. An illustration 
    of this setting can be found in Figure \ref{fig:Theorem4MilnorTubes}.

    \begin{figure}[h]
        \centering
        \includegraphics[scale=0.8]{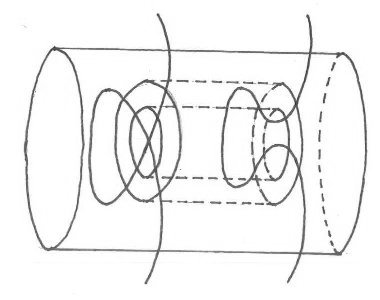}
        \caption{$\overline Y$ and the Milnor tubes $S$ and $T$}
        \label{fig:Theorem4MilnorTubes}
    \end{figure}

    By the Ehresmann 
    fibration theorem we can choose a differentiable flow
    \[
        \Phi: D \times U' \to U' 
    \] 
    defined on an open set $0\in D \subset \mathbb C$ in the deformation base
    and a neighborhood $U'$ of the closure $\overline U$ in 
    $\overline Y$, such that the restriction to $\overline U_0$ 
    \[ 
        \Phi|_{\{\varepsilon\} \times \overline U_0} : 
        \overline U_0 \to \overline U_\varepsilon 
    \] 
    and the restriction to the overlap
    \[
        \Phi|_{\{\varepsilon\}\times \overline W_0} : 
        \overline W_0 \to \overline W_\varepsilon.
    \]
    are diffeomorphisms of manifolds with boundary for sufficiently small
    $|\varepsilon|$. 
    Now consider the long exact sequence in reduced homology 
    for the pair of spaces
    $(\overline Y_0, \overline V_0)$:
    \begin{align}
        \begin{xy}
            \xymatrix{
                \dots \ar[r] &
                H_{i}(\overline V_0) \ar[r] &
                H_{i}(\overline Y_0) \ar[r] & 
                H_i( \overline Y_0, \overline V_0) \ar[r] & 
                H_{i-1}(\overline V_0 ) \ar[r] & 
                \dots
            }
        \end{xy}
        \label{eqn:LESpairOfSpaces1}
    \end{align}
    Because $\overline V_0$ is contractible to a union of points,
    we find 
    $H_i(\overline Y_0) \to H_i(\overline Y_0, \overline V_0)$ 
    to be an isomorhpism for $i > 0$. 
    Clearly, 
    \begin{align}
        H_i(\overline Y_0,\overline V_0 ) = H_i( \overline Y_0/ \overline V_0) 
        = H_i (\overline U_0 / \partial \overline U_0) 
        = H_i (\overline U_\varepsilon / \partial \overline U_\varepsilon )
        = H_i (\overline Y_\varepsilon, \overline V_\varepsilon ),
        \label{eqn:IsomRelativeHomologyGroups}
    \end{align}
    because of excision.
    The identifications in equation \ref{eqn:IsomRelativeHomologyGroups} and 
    equation \ref{eqn:BettiNumbersTjurinaTransfSimpleSing} provide
    us with zeros for the terms 
    $H_i(\overline Y_\varepsilon, \overline V_\varepsilon) = 0$ for 
    $i \neq 0,2$.
    Now consider the analogous long exact sequence for the pair of spaces 
    $(\overline Y_\varepsilon, \overline V_\varepsilon)$ and obtain in degree 
    $3$ 
    \begin{equation}
        \begin{xy}
            \xymatrix{
                0 \ar[r] &
                H_3(\overline V_\varepsilon) \ar[r] &
                H_3(\overline Y_\varepsilon) \ar[r] &
                0
            },
        \end{xy}
        \label{eqn:IdentifDegreeThree}
    \end{equation}
    which means that every vanishing cycle of the occurring singularities
    in $Y_0$ is preserved in the whole fiber $X_\varepsilon$.
    Since for ICIS singularities the Milnor fiber has the topological type
    of a bouquet of $\mu$ $3$-spheres and we have a decomposition
    \[ H_3(\overline V_\varepsilon) 
        = \bigoplus_{i=1}^n H_3(\overline V_{i,\varepsilon}),
    \]
    the middle Betti number $b_3$ of $\overline Y_\varepsilon$ is 
    the sum of Milnor numbers of the singularities of $Y_0$.

    The Tjurina transform $Y_0$ has only ICIS, so we 
    get zeros for $H_i(\overline V_\varepsilon)$ for $i \neq 0,3$ 
    which leads to
    \begin{equation}
        \begin{xy}
            \xymatrix{
                0 \ar[r] &
                H_2(\overline Y_\varepsilon) \ar[r] & 
                H_2(\overline Y_\varepsilon, \overline V_\varepsilon ) \ar[r] &
                0 
            },
        \end{xy}
        \label{eqn:IdentifDegreeTwo}
    \end{equation}
    in degree $2$. 
    As the fibers $\overline Y_\varepsilon$ and $\overline X_\varepsilon$ are 
    isomorphic according to Proposition 
    \ref{prp:IsomorphismSmoothFibersTjurinaModification}, the claim is proved.
\end{proof}

\begin{example}
    The applied ideas also work for higher Cohen-Macaulay type, as we 
    would like to illustrate by revisiting our previous Example 
    \ref{exp:3by4MatrixInFamily} from Section 3: 
    also in the case $t=3$ Proposition 
    \ref{prp:FlatnessTjurinaModificationSimpleSingularities}
    assures the Tjurina modification to work in family.
    Contrary to the case $t=2$, the central fiber in the Tjurina
    modification $Y_0$ now has the homotopy type of 
    $E = \pi_0^{-1}(\{0\}) = \mathbb P^2$. Thus the Betti numbers read
    \begin{equation}
        b_0(Y_0) = 1, \quad
        b_1(Y_0) = 0, \quad
        b_2(Y_0) = 1, \quad
        b_3(Y_0) = 0, \quad
        b_4(Y_0) = 1.
        \label{eqn:BettiNumbers3by4Example}
    \end{equation}
    Clearly the $4$-cycle generated by the $\mathbb P^2$ itself
    can not be preserved in a smooth fiber $Y_\varepsilon \cong X_\varepsilon$,
    because $X_\varepsilon$ is affine and we would get 
    a contradiction to the Lefschetz Hyperplane Theorem.
    In fact this cycle breaks at the 10 points of 
    the $A_1$ singularities of $Y_0$ when resolving them. 
    We can again calculate along the lines of the proof of 
    Theorem \ref{thm:TjTransfInterplay} using the long exact
    sequence of pairs of spaces. 
    The degree $2$ part stays isolated, so we again get
    \[
        b_2(X_\varepsilon) = 1.
    \]
    For degree $3$ the calculations show, that the 
    broken $\mathbb P^2$ leads to a relation among the 
    vanishing cycles of the $A_1$'s. Hence we have
    \[
        b_3(X_\varepsilon) = 9
    \]
    and not $10$ as one might have expected.
\end{example}

\begin{remark}
    A similar phenomenon can be observed when projectivizing 
    the column space of the matrix $M$ of an ICMC2 threefold $(X_0,0)$ of 
    Cohen-Macaulay type $t=2$. As in the case of a 
    Tjurina modification, projectivizing the column space is compatible
    with deformations, because we again get a locally 
    complete intersection transform, say $Z_0$. Only this time
    we find a $\mathbb P^2$ as exceptional set. 
    If we apply this to the first entry of the table 
    \ref{tab:ResultsThreefolds}, also called the $A_0^+$ 
    singularity, we find 
    \[
        b_0(Z_0) = 1, \quad
        b_1(Z_0) = 0, \quad
        b_2(Z_0) = 1, \quad
        b_3(Z_0) = 0, \quad
        b_4(Z_0) = 1
    \]
    and one $A_1$ singularity in $Z_0$. 
    Again, the induced deformation destroys the $4$-cycle of $Z_0$
    leading to the vanishing cycle of the $A_1$ singularity
    being homologous to zero in $X_\varepsilon = Z_\varepsilon$.
\end{remark}

%% file: chapter5.tex
\section{The topological type of the simple ICMC2 singularities}

Using the results of the previous section, direct computation now provides
explicit results for the structure of the Milnor fiber for the simple ICMC2
singularities and for the bounding non-simple ones. We have summarized the
results in the following two tables. Subsequently, we finish this article
by pointing out and explaining some notable observations and stating some
arising questions.

\input{table_results.tex}
\input{table_bounding_singularities.tex}

\begin{remark}{(direct observations from the table)}
\begin{enumerate}
\item We only see simple singularities
occurring in $Y_0$ in table (\ref{tab:ResultsThreefolds}). In table 
(\ref{tab:boundingThreefolds}), where the listed singularities are non-simple, 
there are some simple and some non-simple singularities arising from 
Tjurina modification. In particular, the non-simple ones arise in the 
cases with 1-jet types $J^{(5,2)}$ and in some subcases of $J^{(4,4)}$ in 
the notation of \cite{FN}, whereas the simple ones occur for $J^{(4,2)}$ and 
the remaining subcases fo $J^{(4,4)}$.\\
The 1-jet types $J^{(4,5)}$ and $J^{(4,6)}$ have not been included in the
table, because they lead to non-isolated singularities in $Y_0$, the singular
locus being the whole exceptional ${\mathbb P}^1$.

\item Looking at the preceding tables, one fact immediately attracts attention:
    For any of the explicitly computed ICMC2 threefold singularities
    the second Betti number is always $1$.
    The mechanism behind this fact can be explained as follows:
    According to the identifications 
    (\ref{eqn:IsomRelativeHomologyGroups}) and 
    (\ref{eqn:IdentifDegreeTwo}),
    the second homology group of the smooth fiber 
    is inherited from the exceptional 
    set in the Tjurina transform, which is a $\PP^{(t-1)}$ depending 
    only on the Cohen-Macaulay type $t=2$. 
    In particular this vanishing cycle can not be directly related 
    to any deformation parameters in the sense of e.g. the 
    L\^e-Greuel formulas. 
\end{enumerate}
\end{remark}

    Our results answer negatively to a question of J. Damon and B. Pike who
    expected both Betti numbers $b_2$ and $b_3$ to grow at the same rate
    in families of ICMC2 with 
    \[ \mu = b_3 - b_2 \]
    constant. But our computations show that there are infinite families in
    which neither of the two changes. The first occurrence 
    of this phenomenon can be found 
    in line 7 of the preceding table (\ref{tab:ResultsThreefolds}), 
    given by the matrix
    \[
        \begin{pmatrix}
            w & y & x \\
            z & w & y + v^k
        \end{pmatrix}
    \]
    Following the notation of \cite{FN}, we call these singularities 
    the $\Pi_k$ family.
    The topological type stays the same within this family, whereas
    the isomorphism classes of the space germs do not coincide for different 
    $k$.\\

    The Tjurina transform $Y_0$ of the $\Pi_k$ singularities is smooth. 
    The second summand of the splitting (\ref{eqn:FirstOrderDefY0}) 
    of the space of first order deformations of $Y_0$ vanishes and 
    we obtain
    \[
        N' = H^1(Y_0,T_{Y_0}).
    \]
    More general if we allow isolated singularities in the Tjurina 
    transform as in Lemma (\ref{lem:ISinTjurinaTransfToJetOfMatrix}), 
    the dimension $h^1(Y_0,T_{Y_0})$ measures how many 
    degrees of freedom for deformations of $(X_0,0)$ become locally 
    trivial for all the singular points in $Y_0$. 

    For all the ICIS at the points in the singular locus $\Sigma(Y_0)$
    of the Tjurina transform $Y_0$, 
    we get local vanishing cycles in degree $3$ 
    as we pass to the smooth fiber. 
    The computations (\ref{eqn:IdentifDegreeThree}) in the proof of 
    Theorem (\ref{thm:TjTransfInterplay})
    show that there 
    can be no global relations between them and that they generate 
    the third homology group of the global Milnor fiber $Y_\varepsilon$. 

    Contrary to the second homology group, in degree three we can use the 
    L\^e-Greuel formulas to relate the third Betti number
    \[
        b_3 = \sum_{p \in \Sigma(Y_0)} \mu_p
    \]
    to those of the 
    local singularities $(Y_0,p)$ in the Tjurina transform 
    and to their Tjurina numbers 
    $\tau_p = \dim T^1_{Y_0,p}$. 
    In case of hypersurface singularities in $Y_0$, the numbers $\tau_p$ and 
    $\mu_p$ can be computed as the vector space dimensions 
    of the Tjurina algebra and the Milnor algebra respectively. 
    This covers the case of all simple ICMC2 singularities in 
    $({\mathbb C}^5,0)$, because due to Corollary 
    (\ref{cor:VersalityDefTjurinaTransf}) we can only find simple 3-dimensional
    ICIS in $Y_0$; according to Giusti's list of simple ICIS, these can only be
    A-D-E singularities, see \cite{Giu}. Since the 
    Milnor and the Tjurina numbers coincide for A-D-E singularities (and
    more generally for quasihomogeneous hypersurface singularities), we 
    have proved the following theorem.

    \begin{theorem}
        For ICMC2 threefold singularities of type $2$, whose 
        Tjurina transform $Y_0$ has as singular locus $\Sigma(Y_0)$ 
        either the empty set or a set of dimension $0$, we have 
        \begin{eqnarray}
            \tau &=&  h^1(Y_0,T_{Y_0}) + \sum_{p\in \Sigma(Y_0)} \tau_p, 
            \label{eqn:DifferenceTauGlobalTauLokal}
        \end{eqnarray}
        which becomes 
        \begin{eqnarray}
            \tau &=& h^1(Y_0,T_{Y_0}) + b_3.
            \label{eqn:DifferenceTauB3}
        \end{eqnarray}
        for the simple ICMC2 singularities.
        \label{thm:TauAndB3}
    \end{theorem}

    There seems to be a relation between the number 
    $h^1(Y_0,T_{Y_0})$ and the maximal number of 
    matrix singularities to which $(X_0,0)$ can be deformed. 
    The final object in the adjacencies among ICMC2 3-fold singularities
    of Cohen-Macaulay type 2 is the so called 
    $A_0^+$ singularity
    \[
        \begin{pmatrix}
            x & y & z \\
            v & w & x
        \end{pmatrix},
    \]
    which is the first entry of the table 
    (\ref{tab:ResultsThreefolds}).
    Consider the $\Pi_k$-family from table \ref{tab:ResultsThreefolds} 
    and the deformation 
    over $\CC$ given by
    \[
        \begin{pmatrix}
            w & y & x \\
            z & w & y + v^k
        \end{pmatrix}
        -
        \begin{pmatrix}
            0 & 0 & 0\\
            0 & 0 & \varepsilon
        \end{pmatrix}.
    \]
    For $\varepsilon \neq 0$ we find $k$ distinct singularities
    at the points 
    $(x,y,z,v,w) = (0,0,0, \sqrt[k]{\varepsilon}, 0)$.
    Using the analytic coordinate change 
    $v' = y + v^k -\varepsilon$ locally at any of these points 
    gives the standard form of the $A_0^+$. 
    This and other examples yield $h^1(Y_0,T_{Y_0})$ to grow 
    linearly with the maximal number of $A_0^+$ singularities 
    on a neighboring singular fiber, an observation which coincides
    with the fact that the Tjurina transform is blind to components of the
    discriminant above which only determinantal singularities exist. This 
    train of thought is pursued in detail in \cite{FK4}.

    \begin{figure}[h]
        \centering
        \includegraphics[scale=0.7]{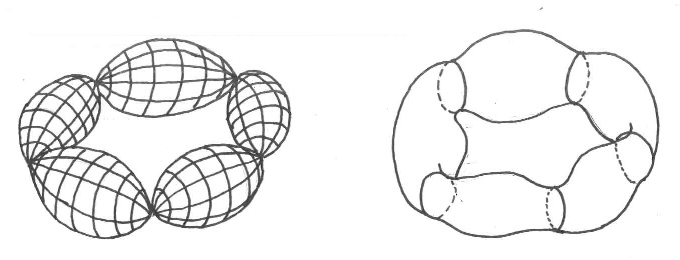}
        \caption{The wrapped candy chain for $k=5$ and $\varepsilon \neq 0$ 
        before and after deforming the $A_0^+$ singularities}
        \label{fig:CandyChain}
    \end{figure}

    The preceding example also exhibits another interesting behaviour in the 
    topology. While the Milnor fiber has its third homology group 
    equal to zero, the fiber over $\varepsilon \neq 0$ contains
    a bouquet of $k$ real $3$-spheres as indicated 
    in Figure \ref{fig:CandyChain}. It reminds us of 
    a chain of wrapped candy.
    It turns out that their total sum is homologous to zero, 
    while either $k-1$ of them generate $H_3(X_\varepsilon)$.
    The $A_0^+$ singularities sit at the ``wrapping points'' 
    between these spheres and are unraveled when passing 
    to a smooth fiber. Consequently 
    all local $2$-cycles become pairwise homologous. 
    A direct computation then shows that the appearing $3$-chain 
    is in fact homologous to $0$.

%% file: table_results.tex
\begin{center}
\begin{longtable}[h]{|c|c|p{2cm}|c|c|}
    \hline
    $M^T$ & $\tau$ & sing. in $Y_0$ & $b_2$ & $b_3$  \\
    \hline
    $
    \begin{pmatrix}
        x & y & z \\
        v & w & x
    \end{pmatrix}
    $ &
    $1$ &
    - &
    $1$ & $0$ \\
    \hline

    $
    \begin{pmatrix}
        x & y & z \\
        v & w & x^{k+1} + y^2
    \end{pmatrix}
    $ & 
    $k+2$ & 
    $A_k$ & 
    $1$ & $k$  \\
    $
    \begin{pmatrix}
        x & y & z \\
        v & w & x y^2 + x^{k-1} 
    \end{pmatrix}
    $ & 
    $k+2$ & 
    $D_{k}$ &
    $1$ & $k$\\

    $
    \begin{pmatrix}
        x & y & z \\
        v & w & x^3 + y^4
    \end{pmatrix}
    $ & 
    $8$ & 
    $E_6$ & $1$ & $6$ \\

    $
    \begin{pmatrix}
        x & y & z \\
        v & w & x^3+xy^3
    \end{pmatrix}
    $ & 
    $9$ & 
    $E_7$ &
    $1$ & $7$ \\

    $
    \begin{pmatrix}
        x & y & z \\
        v & w & x^3 + y^5
    \end{pmatrix}
    $ & 
    $10$ & 
    $E_8$ &
    $1$ & $8$ \\
    \hline

    $
    \begin{pmatrix}
        w & y & x \\
        z & w & y+v^k
    \end{pmatrix}
    $ &
    $2k-1$ & 
    - & 
    $1$ & $0$ \\

    $
    \begin{pmatrix}
        w & y & x \\
        z & w & y^k + v^2
    \end{pmatrix}
    $ &
    $k+2$ & 
    $A_{k-1}$ &
    $1$ & $k-1$ \\

    $
    \begin{pmatrix}
        w & y & x \\
        z & w & yv + v^k 
    \end{pmatrix}
    $ &
    $2k$ & 
    $A_1$ &
    $1$ & $1$ \\

    $
    \begin{pmatrix}
        w+v^k & y & x \\
        z & w & yv
    \end{pmatrix} 
    $ &
    $2k+1$ &
    $A_1$ &
    $1$ & $1$ \\

    $
    \begin{pmatrix}
        w+v^2 & y & x \\
        z & w & y^2 + v^k
    \end{pmatrix}
    $ & 
    $k+3$ & 
    $A_{k-1}$ &
    $1$ & $k-1$ \\

    $
    \begin{pmatrix}
        w & y & x \\
        z & w & y^2 + v^3
    \end{pmatrix}
    $ &
    $7$ &
    $A_2$ & 
    $1$ & $2$\\

    \hline

    $
    \begin{pmatrix}
        v^2 + w^k & y & x \\
        z & w & v^2 + y^l
    \end{pmatrix}
    $
    & $k+l+1$ &
    $A_{k-1}$, $A_{l-1}$ & 
    $1$ & $k+l-2$ \\

    $
    \begin{pmatrix}
        v^2 + w^k & y & x \\
        z & w & yv
    \end{pmatrix}
    $
    & $k+4$ &
    $A_{k-1}$, $A_1$ &
    $1$ & $k$ \\

    $
    \begin{pmatrix}
        v^2 + w^k & y & x \\
        z & w & y^2+v^l
    \end{pmatrix}
    $
    & $k+l+2$ &
    $A_{k-1}$, $A_{l-1}$ &
    $1$ & $k+l-2$ \\ 

    \hline
    $
    \begin{pmatrix}
        wv+ v^k & y & x \\
        z & w & yv + v^k
    \end{pmatrix}
    $
    & $2k+1$ &
    $A_1$, $A_1$ & 
    $1$ & $2$ \\

    $
    \begin{pmatrix}
        wv + v^k & y & x \\
        z & w & yv
    \end{pmatrix}
    $
    & $2k+2$ &
    $A_1$, $A_1$ & 
    $1$ & $2$ \\
    \hline

    $
    \begin{pmatrix}
        wv + v^3 & y & x \\
        z & w & y^2+v^3
    \end{pmatrix}
    $
    & $8$ &
    $A_1$, $A_2$ & 
    $1$ & $3$ \\

    $
    \begin{pmatrix}
        wv & y & x \\
        z & w & y^2 + v^3
    \end{pmatrix}
    $
    & $9$ &
    $A_1$, $A_2$ &
    $1$ & $3$ \\

    $
    \begin{pmatrix}
        w^2 + v^3 & y & x \\
        z & w & y^2+v^3
    \end{pmatrix}
    $
    & $9$ &
    $A_2$, $A_2$ &
    $1$ & $4$ \\
    \hline

    $
    \begin{pmatrix}
        z & y & x \\
        x & w & v^2 + y^2 + z^k
    \end{pmatrix}
    $ 
    & 
    $k+4$ &
    $D_{k+1}$ &
    $1$ & $k+1$ \\

    $
    \begin{pmatrix}
        z & y & x \\
        x & w & v^2 + yz + y^k w
    \end{pmatrix}
    $ & 
    $2k+5$ & 
    $A_{2k+2}$ &
    $1$ & $2k+2$ \\

    $
    \begin{pmatrix}
        z & y & x \\
        x & w & v^2 + yz + y^{k+1}
    \end{pmatrix}
    $ & 
    $2k+4$ & 
    $A_{2k+1}$  &
    $1$ & $2k+1$ \\

    $
    \begin{pmatrix}
        z & y & x \\
        x & w & v^2 + yw + z^2
    \end{pmatrix}
    $ & 
    $8$ &
    $D_5$ & 
    $1$ & $5$ \\

    $
    \begin{pmatrix}
        z & y & x \\
        x & w & v^2 + y^3 + z^2
    \end{pmatrix}
    $ & 
    $9$ &
    $E_6$ &
    $1$ & $6$ \\

    $
    \begin{pmatrix}
        z & y & x + v^2 \\
        x & w & vy + z^2
    \end{pmatrix}
    $ & 
    $7$ &
    $D_3$ &
    $1$ & $3$ \\

    $
    \begin{pmatrix}
        z & y & x + v^2 \\
        x & w & vz + y^2
    \end{pmatrix}
    $ & 
    $8$ & 
    $A_4$ &
    $1$ & $4$ \\

    $
    \begin{pmatrix}
        z & y & x + v^2 \\
        x & w & z^2 + y^2
    \end{pmatrix}
    $ & 
    $9$ &
    $D_5$ &
    $1$ & $5$ \\
    \hline
    \caption{Homology of Milnor fibers computed by means of the Tjurina 
             modification}
    \label{tab:ResultsThreefolds}
\end{longtable}
    
\end{center}

%% file: table_bounding_singularities.tex
\begin{center}
\begin{longtable}[h]{|c|c|p{2cm}|c|c|}
    \hline
    $M^T$ & $\tau$ & sing. in $Y_0$ & $b_2$ & $b_3$ \\
    \hline
    $
    \begin{pmatrix}
        x & y & z \\
        w & v & x^4 + y^4 
    \end{pmatrix}
    $ &
    $11$ &
    $X_9$ &
    $1$ & $9$ \\
    \hline

    $
    \begin{pmatrix}
        x & y & z \\
        w & v & x^3 + y^6
    \end{pmatrix}
    $ &
    $12$ &
    $J_{10}$ &
    $1$ & $10$ \\
    \hline

    $
    \begin{pmatrix}
        w+v^2 & y & x \\
        z & v & y^3 + v^3

    \end{pmatrix}
    $ &
    $8$ &
    $D_4$ &
    $1$ & $4$ \\
    \hline

    $
    \begin{pmatrix}
        w+v^3 & y & x \\
        z & w & y^2 + v^4
    \end{pmatrix}
    $ &
    $9$ &
    $A_3$ &
    $1$ & $3$ \\
    \hline

    $
    \begin{pmatrix}
        z & y & x \\
        x & w & v^2 + y^3 + z^3
    \end{pmatrix}
    $ &
    $11$ &
    $T_{3,3,3}$ &
    $1$ & $8$ \\
    \hline

    $
    \begin{pmatrix}
        z & y & x \\
        x & w & v^3 + y^2 + z^3
    \end{pmatrix}
    $\footnote{There is a typesetting error in this matrix in \cite{FN}. The
               right-hand lower entry here is the correct one.} &
    $13$ &
    $T_{3,3,3}$ &
    $1$ & $8$ \\
    \hline

    $
    \begin{pmatrix}
        z & y & x\\
        x & w & v^3 + y^3 + z^2
    \end{pmatrix}
    $ &
    $17$ &
    $U_{12}$ &
    $1$ & $12$ \\
    \hline

    $
    \begin{pmatrix}
        z & y & x \\
        x & w & v^2 + y^4 + z^2
    \end{pmatrix}
    $ &
    $12$ &
    $X_9$ &
    $1$ & $9$ \\
    \hline

    $
    \begin{pmatrix}
        z & y & x+v^2 \\
        x & w & vz + yz + vw
    \end{pmatrix}
    $ &
    $10$ &
    $D_6$ &
    $1$ & $6$ \\
    \hline

    $
    \begin{pmatrix}
        z & y & x+v^3 \\
        x & w & vy+z^2
    \end{pmatrix}
    $ &
    $9$ &
    $A_3$ &
    $1$ & 3 \\
    \hline

    $
    \begin{pmatrix}
        z & y & x+v^3 \\
        x & w & y^2 + yz + z^2
    \end{pmatrix}
    $ &
    $15$ &
    $X_9$ &
    $1$ & $9$ \\
    \hline

    $
    \begin{pmatrix}
        z & y & x+v^2 \\
        x & w & vy + yz + z^3
    \end{pmatrix}
    $ &
    $8$ &
    $D_4$ &
    $1$ & 4 \\
    \hline
    \caption{Homology of Milnor fibers for the bounding non-simple singularities}
    \label{tab:boundingThreefolds}
\end{longtable}
    
\end{center}